\documentclass[11pt]{amsart}
\def\N{{\mathbb{N}}}
\def\Z{{\mathbb{Z}}}
\def\R{{\mathbb{R}}}
\def\Image{\operatorname{Im}}
\def\Int{\operatorname{Int}}
\def\diam{\operatorname{diam}}
\newtheorem{Theorem}{Theorem}[section]
\newtheorem{Corollary}[Theorem]{Corollary}
\newtheorem{Lemma}[Theorem]{Lemma}

\newtheorem{Thm}{Theorem}
\newtheorem{Cor}[Thm]{Corollary}
\theoremstyle{definition}
\newtheorem{Definition}[Theorem]{Definition}
\newtheorem{Example}[Theorem]{Example}

\theoremstyle{remark}

\begin{document}
\sloppy
\title{CAT(0) groups and Coxeter groups whose boundaries are scrambled sets}
\author{Tetsuya Hosaka} 
\address{Department of Mathematics, Faculty of Education, 
Utsunomiya University, Utsunomiya, 321-8505, Japan}
\date{April 28, 2007}
\email{hosaka@cc.utsunomiya-u.ac.jp}
\thanks{
Partly supported by the Grant-in-Aid for Young Scientists (B), 
The Ministry of Education, Culture, Sports, Science and Technology, Japan.
(No.\ 18740025).}
\maketitle
\begin{abstract}
In this paper, 
we study CAT(0) groups and Coxeter groups whose 
boundaries are scrambled sets.
Suppose that a group $G$ acts 
geometrically (i.e.\ properly and cocompactly by isometries) 
on a CAT(0) space $X$.
(Such group $G$ is called a {\it CAT(0) group}.)
Then the group $G$ acts by homeomorphisms 
on the boundary $\partial X$ of $X$ and 
we can define a metric $d_{\partial X}$ on the boundary $\partial X$.
The boundary $\partial X$ is called a {\it scrambled set} 
if for any $\alpha,\beta\in\partial X$ with $\alpha\neq\beta$, 
(1) $\limsup\{d_{\partial X}(g\alpha,g\beta)\,|\,g\in G\}>0$ and 
(2) $\liminf\{d_{\partial X}(g\alpha,g\beta)\,|\,g\in G\}=0$.
We investigate when are boundaries of CAT(0) groups 
(and Coxeter groups) scrambled sets.
\end{abstract}

\section{Introduction}

The purpose of this paper is to study 
CAT(0) groups and Coxeter groups whose 
boundaries are scrambled sets.

Definitions and basic properties of CAT(0) spaces and their boundaries 
are found in \cite{BH}.
A {\it geometric} action on a CAT(0) space 
is an action by isometries which is proper (\cite[p.131]{BH}) and cocompact.
We note that every CAT(0) space 
on which some group acts 
geometrically is a proper space (\cite[p.132]{BH}).
A group $G$ is called a {\it CAT(0) group}, 
if there exists a geometric action of $G$ on some CAT(0) space $X$.
Here we say that 
the boundary $\partial X$ of $X$ is a {\it boundary of $G$}.
We note that 
if $G$ is hyperbolic then the group $G$ determines 
the boundary $\partial X$.
In general, a CAT(0) group $G$ does not determine 
the boundary $\partial X$ of a CAT(0) space $X$ 
on which $G$ acts geometrically (Croke and Kleiner \cite{CK}).

Suppose that a group $G$ acts geometrically on a CAT(0) space $X$.
Then the group $G$ acts by homeomorphisms 
on the boundary $\partial X$ of $X$.
We note that if $|\partial X|>2$ then 
the boundary $\partial X$ is uncountable.

We define a metric on the boundary $\partial X$ as follows: 
We first fix a basepoint $x_0\in X$.
Let $\alpha,\beta\in\partial X$.
There exist unique geodesic rays 
$\xi_{x_0,\alpha}$ and $\xi_{x_0,\beta}$ in $X$ 
with $\xi_{x_0,\alpha}(0)=\xi_{x_0,\beta}(0)=x_0$, 
$\xi_{x_0,\alpha}(\infty)=\alpha$ and $\xi_{x_0,\beta}(\infty)=\beta$.
Then the metric $d_{\partial X}(\alpha,\beta)$ 
of $\alpha$ and $\beta$ on $\partial X$ 
(with respect to the basepoint $x_0$)
is defined by 
$$ d_{\partial X}(\alpha,\beta)
=\sum_{i=1}^{\infty}\min\{d(\xi_{x_0,\alpha}(i),\xi_{x_0,\beta}(i)),\ 
\frac{1}{2^i}\}.$$
The metric $d_{\partial X}$ depends on the basepoint $x_0$ and 
the topology of $\partial X$ does not depend on $x_0$.

The boundary $\partial X$ is said to be {\it minimal}, 
if any orbit $G\alpha$ is dense in $\partial X$.
Also the boundary $\partial X$ is called a {\it scrambled set}, 
if for any $\alpha,\beta\in\partial X$ with $\alpha\neq\beta$, 
\begin{align*}
&\limsup\{d_{\partial X}(g\alpha,g\beta)\,|\,g\in G\}>0 \ \text{and}\\
&\liminf\{d_{\partial X}(g\alpha,g\beta)\,|\,g\in G\}=0.
\end{align*}

We note that ``minimality'' and ``scrambled sets'' in this paper 
is a natural extension of the original definitions in the chaotic theory 
and original minimality and scrambled sets 
are important concept on dynamical systems and 
relate to the chaotic theory in the sense of Li and Yorke 
(cf.\ \cite{HuaY}, \cite{Ka} and \cite{LY}).
Boundaries of CAT(0) groups (and Coxeter groups) are 
interesting object.
In general, they are so complex 
and it is so difficult to see these topology and 
the actual actions of CAT(0) groups on their boundaries.
The purpose of this paper is to 
get something information of these actions and boundaries 
by using a method of the chaotic theory.
We can find recent research on minimality 
of boundaries of CAT(0) groups and Coxeter groups 
in \cite{Ho1}, \cite{Ho2}, \cite{Ho3} and \cite{Ho5}.
In this paper, 
we investigate CAT(0) groups and Coxeter groups 
whose boundaries are scrambled sets.

After some preliminaries on CAT(0) spaces and their boundaries in Section~2, 
we first show the following theorem in Section~3.

\begin{Thm}\label{Thm1}
Suppose that a group $G$ acts 
geometrically on a CAT(0) space $X$ and $|\partial X|>2$.
Then 
$$\limsup\{d_{\partial X}(g\alpha,g\beta)\,|\,g\in G\}>0$$ 
for any $\alpha,\beta\in\partial X$ with $\alpha\neq\beta$.
\end{Thm}

This theorem implies that 
the boundary $\partial X$ is a scrambled set if and only if 
$$\liminf\{d_{\partial X}(g\alpha,g\beta)\,|\,g\in G\}=0$$ 
for any $\alpha,\beta\in\partial X$ with $\alpha\neq\beta$.

In Section~3, 
we prove a technical theorem 
which gives a sufficient condition of CAT(0) groups 
whose boundaries are scrambled sets 
and which plays a key role 
in the proof of the main results in this paper.

In Section~4, 
we study boundaries of hyperbolic CAT(0) groups 
and we show the following theorems.

\begin{Thm}\label{Thm3}
The boundary of every non-elementary hyperbolic CAT(0) group 
is a scrambled set.
\end{Thm}

An action of a group $G$ on a metric space $Y$ by homeomorphisms 
is said to be {\it expansive}, 
if there exists a constant $c>0$ such that 
for each pair $y,y'\in Y$ with $y\neq y'$, 
there is $g\in G$ such that $d(gy,gy')>c$.

\begin{Thm}\label{Thm3-2}
Suppose that a group $G$ acts 
geometrically on a CAT(0) space $X$ and $|\partial X|>2$.
The action of $G$ on $\partial X$ is expansive 
if and only if the space $X$ is hyperbolic.
\end{Thm}

In Section~5, 
we investigate when are boundaries of CAT(0) groups scrambled sets, 
and we give a sufficient condition of 
CAT(0) groups whose boundaries are scrambled sets and 
also give a sufficient condition of CAT(0) groups 
whose boundaries are {\it not} scrambed sets.

In Sections~6, 7 and 8, 
we study the boundary of a Coxeter system.
Definitions and basic properties of Coxeter systems and Coxeter groups 
are found in \cite{Bo} and \cite{Hu}.
Every Coxeter system $(W,S)$ determines 
a {\it Davis complex} $\Sigma(W,S)$ 
which is a CAT(0) space (\cite{D1}, \cite{D2}, \cite{D3}, \cite{M}) 
and the natural action of the Coxeter group $W$ 
on the Davis complex $\Sigma(W,S)$ 
is proper, cocompact and by isometries 
(hence Coxeter groups are CAT(0) groups).
The boundary $\partial\Sigma(W,S)$ is called the {\it boundary} of 
the Coxeter system $(W,S)$.

We show a technical theorem 
which gives a sufficient condition of a Coxeter system 
whose boundary is a scrambled set in Section~7.

Using the technical theorem, 
we show the following strong theorem for right-angled Coxeter groups 
and their boundaries in Section~8.

\begin{Thm}\label{Thm6}
If $(W,S)$ is an irreducible right-angled Coxeter system 
and $|\partial\Sigma(W,S)|>2$, 
then the boundary $\partial\Sigma(W,S)$ is a scrambled set.
\end{Thm}

From Theorem~\ref{Thm6} and \cite[Theorem~5.1]{Ho5}, 
we obtain the following corollary 
which gives equivalent conditions of 
a right-angled Coxeter system 
whose boundary is a scrambled set.

\begin{Cor}\label{Cor7}
Let $(W,S)$ be a right-angled Coxeter system with $|\partial\Sigma(W,S)|>2$.
Then the following statements are equivalent:
\begin{enumerate}
\item[(1)] $\partial\Sigma(W,S)$ is a scrambled set.
\item[(2)] $\partial\Sigma(W,S)$ is minimal.
\item[(3)] $(W_{\tilde{S}},\tilde{S})$ is irreducible.
\end{enumerate}
\end{Cor}

Here $W_{\tilde{S}}$ is the minimum parabolic subgroup 
of finite index in $(W,S)$, that is, 
for the irreducible decomposition 
$W=W_{S_1}\times\dots\times W_{S_n}$, 
we define 
$\tilde{S}=\bigcup\{S_i\,|\,\text{$W_{S_i}$ is infinite}\}$ (\cite{Ho0}) 
and $W_{\tilde{S}}$ is the subgroup of $W$ generated by $\tilde{S}$.

By Corollary~\ref{Cor7}, 
we can determine the class of right-angled Coxeter systems 
whose boundaries are scrambled sets.

\section{CAT(0) spaces and their boundaries}

We say that a metric space $(X,d)$ is a {\it geodesic space} if 
for each $x,y \in X$, 
there exists an isometric embedding $\xi:[0,d(x,y)] \rightarrow X$ such that 
$\xi(0)=x$ and $\xi(d(x,y))=y$ (such $\xi$ is called a {\it geodesic}).
Also a metric space $X$ is said to be {\it proper} 
if every closed metric ball is compact.

Let $X$ be a geodesic space and 
let $T$ be a geodesic triangle in $X$.
A {\it comparison triangle} for $T$ is 
a geodesic triangle $\overline{T}$ in the Euclidean plane $\R^2$
with same edge lengths as $T$.
Choose two points $x$ and $y$ in $T$. 
Let $\bar{x}$ and $\bar{y}$ denote 
the corresponding points in $\overline{T}$.
Then the inequality $$d(x,y) \le d_{\R^2}(\bar{x},\bar{y})$$ 
is called the {\it CAT(0)-inequality}, 
where $d_{\R^2}$ is the usual metric on $\R^2$.
A geodesic space $X$ is called a {\it CAT(0) space} 
if the CAT(0)-inequality holds
for all geodesic triangles $T$ and for all choices of two points $x$ and 
$y$ in $T$.

Let $X$ be a proper CAT(0) space and $x_0 \in X$.
The {\it boundary of $X$ with respect to $x_0$}, 
denoted by $\partial_{x_0}X$, is defined as 
the set of all geodesic rays issuing from $x_0$. 
Then we define a topology on $X \cup \partial_{x_0}X$ 
by the following conditions: 
\begin{enumerate}
\item[(1)] $X$ is an open subspace of $X \cup \partial_{x_0}X$. 
\item[(2)] For $\alpha \in \partial_{x_0}X$ and $r, \epsilon >0$, let
$$ U_{x_0}(\alpha;r,\epsilon)=
\{ x \in X \cup \partial_{x_0} X \,|\, 
x \not\in B(x_0,r),\ d(\alpha(r),\xi_{x}(r))<\epsilon \}, $$
where $\xi_{x}:[0,d(x_0,x)]\rightarrow X$ is the geodesic from $x_0$ to $x$
($\xi_{x}=x$ if $x \in \partial_{x_0} X$).
Then for each $\epsilon_0>0$, 
the set 
$$\{U_{x_0}(\alpha;r,\epsilon_0)\,|\, r>0\}$$ 
is a neighborhood basis for $\alpha$ in $X \cup \partial_{x_0}X$. 
\end{enumerate}
This is called the {\it cone topology} on $X \cup \partial_{x_0} X$.
It is known that 
$X \cup \partial_{x_0} X$ is 
a metrizable compactification of $X$ (\cite{BH}, \cite{GH}).

Let $X$ be a geodesic space.
Two geodesic rays $\xi, \zeta:[0,\infty) \rightarrow X$ are 
said to be {\it asymptotic} if there exists a constant $N$ such that 
$d(\xi(t),\zeta(t)) \le N$ for any $t \ge 0$. 
It is known that 
for each geodesic ray $\xi$ in $X$ and 
each point $x \in X$, 
there exists a unique geodesic ray $\xi'$ issuing from $x$ 
such that $\xi$ and $\xi'$ are asymptotic.

Let $x_0$ and $x_1$ be two points of a proper CAT(0) space $X$.
Then there exists a unique bijection 
$\Phi:\partial_{x_0}X \rightarrow \partial_{x_1}X$ 
such that $\xi$ and $\Phi(\xi)$ are asymptotic 
for any $\xi \in \partial_{x_0}X$. 
It is known that $\Phi:\partial_{x_0}X\rightarrow \partial_{x_1}X$ 
is a homeomorphism (\cite{BH}, \cite{GH}).

Let $X$ be a proper CAT(0) space.
The asymptotic relation is an equivalence relation 
in the set of all geodesic rays in $X$.
The {\it boundary of} $X$, denoted by $\partial X$, 
is defined as the set of asymptotic equivalence classes of geodesic rays.
The equivalence class of a geodesic ray $\xi$ is denoted by $\xi(\infty)$.
For each $x_0 \in X$ and each $\alpha \in \partial X$, 
there exists a unique element $\xi \in \partial_{x_0}X$ 
with $\xi(\infty)=\alpha$.
Thus we may identify $\partial X$ with $\partial_{x_0}X$ for each $x_0 \in X$.

We can define a metric on the boundary $\partial X$ as follows: 
We first fix a besepoint $x_0\in X$.
Let $\alpha,\beta\in\partial X$ and 
let $\xi_{x_0,\alpha}$ and $\xi_{x_0,\beta}$ be the geodesic rays in $X$ 
with $\xi_{x_0,\alpha}(0)=\xi_{x_0,\beta}(0)=x_0$, 
$\xi_{x_0,\alpha}(\infty)=\alpha$ and $\xi_{x_0,\beta}(\infty)=\beta$.
Then the metric $d_{\partial X}(\alpha,\beta)$ 
of $\alpha$ and $\beta$ on $\partial X$ is defined by 
$$ d_{\partial X}(\alpha,\beta)
=\sum_{i=1}^{\infty}\min\{d(\xi_{x_0,\alpha}(i),\xi_{x_0,\beta}(i)),\ 
\frac{1}{2^i}\}.$$
We note that the metric $d_{\partial X}$ depends on the basepoint $x_0$.

In this paper, we suppose that 
every CAT(0) space $X$ has a fixed basepoint $x_0$ 
and the metric $d_{\partial X}$ on the boundary $\partial X$ 
is defined by the basepoint $x_0$.

Let $X$ be a proper CAT(0) space 
and let $G$ be a group which acts on $X$ by isometries.
For each element $g \in G$ and 
each geodesic ray $\xi:[0,\infty)\rightarrow X$, 
a map $g \xi:[0,\infty)\rightarrow X$ 
defined by $(g\xi)(t):=g(\xi(t))$ is also a geodesic ray.
If geodesic rays $\xi$ and $\xi'$ are asymptotic, 
then $g\xi$ and $g\xi'$ are also asymptotic.
Thus $g$ induces a homeomorphism of $\partial X$, 
and the group $G$ acts by homeomorphisms on the boundary $\partial X$.

A {\it geometric} action on a CAT(0) space 
is an action by isometries which is proper (\cite[p.131]{BH}) 
and cocompact.
We note that every CAT(0) space 
on which a group acts 
geometrically is a proper space (\cite[p.132]{BH}).
A group which acts geometrically on some CAT(0) space is 
called a {\it CAT(0) group}. 

Details of CAT(0) spaces and their boundaries are found in 
\cite{BH} and \cite{GH}.

Here we introduce some properties of CAT(0) spaces and their boundaries 
used later.

\begin{Lemma}[\cite{BH}, \cite{GH}]\label{lem00}
Let $X$ be a proper CAT(0) space.
\begin{enumerate}
\item[(1)] 
For each three points $x_0,x_1,x_2\in X$ and each $t \in [0,1]$, 
$$d(\xi_1(td(x_0,x_1)),\xi_2(td(x_0,x_2))) \le t d(x_1,x_2),$$
where $\xi_i:[0,d_i]\rightarrow X$ is 
the geodesic segment from $x_0$ to $x_i$ for each $i=1,2$.
\item[(2)] If geodesic rays 
$\xi,\xi':[0,\infty)\rightarrow X$ are asymptotic, then 
$$d(\xi(t), \xi'(t))\le d(\xi(0),\xi'(0))$$ 
for any $t\ge 0$.
\end{enumerate}
\end{Lemma}

From Lemma~\ref{lem00}~(1), we obtain the following.

\begin{Lemma}[\cite{BH}, \cite{GH}]\label{lem0}
Let $X$ be a CAT(0) space and 
let $\xi:[0,\infty)\rightarrow X$ and $\xi':[0,\infty)\rightarrow X$ be 
two geodesic rays with $\xi(0)=\xi'(0)$.
Then for $0<t\le t'$, 
$$ d(\xi(t),\xi'(t))\le \frac{t}{t'}d(\xi(t'),\xi'(t')).$$
\end{Lemma}

We obtain the following lemma from the proof of \cite[Lemma~4.2]{Ho00}.

\begin{Lemma}[{\cite[Lemma~4.2]{Ho00}}]\label{lem1}
Let $X$ be a CAT(0) space and 
let $\xi:[0,\infty)\rightarrow X$ and $\xi':[0,\infty)\rightarrow X$ be 
two geodesic rays with $\xi(0)=\xi'(0)$.
For $r>\epsilon>0$, 
if $d(\xi(r),\Image\xi')\le \epsilon$ then 
$d(\xi(r-\epsilon),\xi'(r-\epsilon))\le \epsilon$.
\end{Lemma}

We define the {\it angle} of two geodesic paths in a CAT(0) space.

\begin{Definition}[{\cite[p.9 and p.184]{BH}}]
Let $X$ be a CAT(0) space and 
let $\xi:[0,a]\rightarrow X$ and $\xi':[0,a']\rightarrow X$ be 
two geodesic paths with $\xi(0)=\xi'(0)$.
For $t\in(0,a]$ and $t'\in(0,a']$, 
we consider the comparison triangle 
$\overline{\triangle}(\xi(0),\xi(t),\xi'(t'))$, and 
the comparison angle $\overline{\angle}_{\xi(0)}(\xi(t),\xi'(t))$.
The {\it Alexandrov angle} between the geodesic paths 
$\xi$ and $\xi'$ is the number $\angle(\xi,\xi')\in [0,\pi]$ defined by 
$$ \angle(\xi,\xi')=\limsup_{t,t'\rightarrow 0}
\overline{\triangle}(\xi(0),\xi(t),\xi'(t')).$$
\end{Definition}

\begin{Lemma}[{\cite[p.184, Proposition~II.3.1]{BH}}]\label{lem:angle0}
Let $X$ be a CAT(0) space and 
let $\xi:[0,a]\rightarrow X$ and $\xi':[0,a']\rightarrow X$ be 
two geodesic paths with $\xi(0)=\xi'(0)$.
Then 
$$ \angle(\xi,\xi')=\lim_{t\rightarrow 0}2\arcsin\frac{1}{2t}d(\xi(t),\xi'(t)).$$
\end{Lemma}

We define the {\it angle} of two points 
in the boundary of a proper CAT(0) space.

\begin{Definition}[{\cite[p.280]{BH}}]
Let $X$ be a proper CAT(0) space, 
let $x\in X$ and let $\alpha,\beta\in\partial X$.
The angle $\angle_x(\alpha,\beta)$ at $x$ between $\alpha$ and $\beta$ 
is defined as 
$$ \angle_x(\alpha,\beta) =\angle_x(\xi_{x,\alpha},\xi_{x,\beta}), $$
where $\xi_{x,\alpha}$ and $\xi_{x,\beta}$ are the geodesic rays 
with $\xi_{x,\alpha}(0)=\xi_{x,\beta}(0)=x$, 
$\xi_{x,\alpha}(\infty)=\alpha$ and $\xi_{x,\beta}(\infty)=\beta$.

Also the angle $\angle(\alpha,\beta)$ between $\alpha$ and $\beta$ 
is defined as 
$$ \angle(\alpha,\beta)=\sup_{x\in X}\angle_x(\alpha,\beta).$$
\end{Definition}

\begin{Lemma}[{\cite[p.281, Proposition~II.9.8]{BH}}]\label{lem:angle}
Let $X$ be a proper CAT(0) space, 
let $x_0\in X$, let $\alpha,\beta\in\partial X$ and 
let $\xi_{x_0,\alpha}$ and $\xi_{x_0,\beta}$ be the geodesic rays 
with $\xi_{x_0,\alpha}(0)=\xi_{x_0,\beta}(0)=x_0$, 
$\xi_{x_0,\alpha}(\infty)=\alpha$ and $\xi_{x_0,\beta}(\infty)=\beta$.
\begin{enumerate}
\item[(1)] 
The function $t\mapsto \angle_{\xi_{x_0,\alpha}(t)}(\alpha,\beta)$ 
is non-decreasing and 
$$ \angle(\alpha,\beta)
=\lim_{t\rightarrow\infty}\angle_{\xi_{x_0,\alpha}(t)}(\alpha,\beta). $$
\item[(2)] 
$\displaystyle 2\sin(\angle(\alpha,\beta)/2)
=\lim_{t\rightarrow\infty}\frac{1}{t}d(\xi_{x_0,\alpha}(t),\xi_{x_0,\beta}(t))$.
\item[(3)] 
$\displaystyle 2\sin(\angle_{x_0}(\alpha,\beta)/2)
=\lim_{t\rightarrow 0}\frac{1}{t}d(\xi_{x_0,\alpha}(t),\xi_{x_0,\beta}(t))$.
\end{enumerate}
\end{Lemma}

Here we obtain (3) in the above lemma from Lemma~\ref{lem:angle0}.

Using Lemmas~\ref{lem0} and \ref{lem:angle}, we show a lemma.

\begin{Lemma}\label{lem:angle1}
Let $X$ be a proper CAT(0) space, 
let $x_0\in X$, let $\alpha,\beta\in\partial X$ and 
let $\xi_{x_0,\alpha}$ and $\xi_{x_0,\beta}$ be the geodesic rays 
with $\xi_{x_0,\alpha}(0)=\xi_{x_0,\beta}(0)=x_0$, 
$\xi_{x_0,\alpha}(\infty)=\alpha$ and $\xi_{x_0,\beta}(\infty)=\beta$.
Then 
$$ 2\sin(\angle_{x_0}(\alpha,\beta)/2) 
\le d(\xi_{x_0,\alpha}(1),\xi_{x_0,\beta}(1)).$$
\end{Lemma}

\begin{proof}
By Lemma~\ref{lem:angle}~(3), 
$$ 2\sin(\angle_{x_0}(\alpha,\beta)/2)
=\lim_{t\rightarrow 0}\frac{1}{t}d(\xi_{x_0,\alpha}(t),\xi_{x_0,\beta}(t)).$$
Here by Lemma~\ref{lem0}, for any $0<t\le t'$ 
$$\frac{1}{t}d(\xi_{x_0,\alpha}(t),\xi_{x_0,\beta}(t)) 
\le \frac{1}{t'}d(\xi_{x_0,\alpha}(t'),\xi_{x_0,\beta}(t')).$$
Hence 
\begin{align*}
2\sin(\angle_{x_0}(\alpha,\beta)/2)
&=\lim_{t\rightarrow 0}\frac{1}{t}d(\xi_{x_0,\alpha}(t),\xi_{x_0,\beta}(t)) \\
&\le d(\xi_{x_0,\alpha}(1),\xi_{x_0,\beta}(1)).
\end{align*}
\end{proof}

\section{A key theorem on CAT(0) groups whose boundaries are scrambled sets}

In this section, 
we show a key theorem which gives a sufficient condition of 
CAT(0) groups whose boundaries are scrambled sets.

We first prove the following theorem.

\begin{Theorem}\label{Thm3-1}
Suppose that a group $G$ acts 
geometrically on a CAT(0) space $X$ and $|\partial X|>2$.
Then 
$$\limsup\{d_{\partial X}(g\alpha,g\beta)\,|\,g\in G\}>0$$ 
for any $\alpha,\beta\in\partial X$ with $\alpha\neq\beta$.
\end{Theorem}

\begin{proof}
Let $\alpha,\beta\in\partial X$ with $\alpha\neq\beta$.
Then $\angle(\alpha,\beta)>0$, since $\alpha\neq\beta$.
Let $\xi_{x_0,\alpha}$ be the geodesic ray 
with $\xi_{x_0,\alpha}(0)=x_0$ and $\xi_{x_0,\alpha}(\infty)=\alpha$ 
and let $x_i=\xi_{x_0,\alpha}(i)$ for each $i\in\N$.
Then by Lemma~\ref{lem:angle}~(1), 
the function $i\mapsto \angle_{x_i}(\alpha,\beta)$ 
is non-decreasing and 
the sequence $\{\angle_{x_i}(\alpha,\beta)\}_i$ 
converges to $\angle(\alpha,\beta)>0$ as $i\rightarrow\infty$.
Hence there exists a number $i_0\in\N$ 
such that for any $i\ge i_0$, 
$$\angle(\alpha,\beta)/2\le\angle_{x_i}(\alpha,\beta)\le\angle(\alpha,\beta).$$
Let $\xi_{x_i,\alpha}$ and $\xi_{x_i,\beta}$ be the geodesic rays 
with $\xi_{x_i,\alpha}(0)=\xi_{x_i,\beta}(0)=x_i$, 
$\xi_{x_i,\alpha}(\infty)=\alpha$ and $\xi_{x_i,\beta}(\infty)=\beta$.
By Lemma~\ref{lem:angle1}, 
we have that 
$$d(\xi_{x_i,\alpha}(1),\xi_{x_i,\beta}(1))\ge 
2\sin(\angle_{x_i}(\alpha,\beta)/2) \ge 2\sin(\angle(\alpha,\beta)/4),$$
because $\angle_{x_i}(\alpha,\beta)\ge\angle(\alpha,\beta)/2$.
Let $r=2\sin(\angle(\alpha,\beta)/4)$.
Then 
\begin{align*}
&d(\xi_{x_i,\alpha}(1),\xi_{x_i,\beta}(1))\ge r\ \text{and} \\
&d(\xi_{x_i,\alpha}(t),\xi_{x_i,\beta}(t))\ge rt
\end{align*}
for any $t\ge 1$ by Lemma~\ref{lem0}.
Since the action of $G$ on $X$ is cocompact and $X$ is proper, 
$GB(x_0,N)=X$ for some $N>0$.
For each $i\in\N$, there exists $g_i\in G$ such that $d(x_i,g_ix_0)\le N$.
Let $\xi_{g_ix_0,\alpha}$ and $\xi_{g_ix_0,\beta}$ be the geodesic rays 
with $\xi_{g_ix_0,\alpha}(0)=\xi_{g_ix_0,\beta}(0)=g_ix_0$, 
$\xi_{g_ix_0,\alpha}(\infty)=\alpha$ and $\xi_{g_ix_0,\beta}(\infty)=\beta$.
By Lemma~\ref{lem00}~(2), 
$d(\xi_{x_i,\alpha}(t),\xi_{g_ix_0,\alpha}(t))\le N$ and 
$d(\xi_{x_i,\beta}(t),\xi_{g_ix_0,\beta}(t))\le N$ for any $t\ge 0$.
Hence 
$$d(\xi_{g_ix_0,\alpha}(t),\xi_{g_ix_0,\beta}(t))\ge
d(\xi_{x_i,\alpha}(t),\xi_{x_i,\beta}(t))-2N \ge rt-2N$$
for each $t\ge 1$.
Let $\displaystyle t_0=[\frac{2N+1}{r}]+1$.
(We note that the number $t_0$ depends on just $\alpha$ and $\beta$.)
Then $rt_0-2N \ge 1$ and 
$$d(\xi_{g_ix_0,\alpha}(t_0),\xi_{g_ix_0,\beta}(t_0))\ge rt_0-2N\ge 1.$$
Here 
$g_i\xi_{x_0,g_i^{-1}\alpha}=\xi_{g_ix_0,\alpha}$ and 
$g_i\xi_{x_0,g_i^{-1}\beta}=\xi_{g_ix_0,\beta}$, 
since $g_i$ is an isometry of $X$.
Hence for each $i\ge i_0$, 
\begin{align*}
d_{\partial X}(g_i^{-1}\alpha,g_i^{-1}\beta)
&=\sum_{j=1}^{\infty}\min\{
d(\xi_{x_0,g_i^{-1}\alpha}(j),\xi_{x_0,g_i^{-1}\beta}(j)),\ 
\frac{1}{2^j}\} \\
&=\sum_{j=1}^{\infty}\min\{
d(g_i^{-1}\xi_{g_ix_0,\alpha}(j),g_i^{-1}\xi_{g_ix_0,\beta}(j)),\ 
\frac{1}{2^j}\} \\
&=\sum_{j=1}^{\infty}\min\{
d(\xi_{g_ix_0,\alpha}(j),\xi_{g_ix_0,\beta}(j)),\ 
\frac{1}{2^j}\} \\
&\ge\frac{1}{2^{t_0}}, 
\end{align*}
because $d(\xi_{g_ix_0,\alpha}(t_0),\xi_{g_ix_0,\beta}(t_0))\ge 1$.
Here $\displaystyle \frac{1}{2^{t_0}}$ is a constant 
which depends on just $\alpha$ and $\beta$.
Thus we obtain that 
$$\limsup\{d_{\partial X}(g\alpha,g\beta)\,|\,g\in G\}
\ge\frac{1}{2^{t_0}}>0.$$
\end{proof}

Now we show a theorem which gives a sufficient condition of CAT(0) groups 
whose boundaries are scrambed sets.
This theorem plays a key role in the proof of the main results in this paper.

\begin{Theorem}\label{Thm:key}
Suppose that a group $G$ acts 
geometrically on a CAT(0) space $X$ and $|\partial X|>2$.
Assume that there exists a constant $M>0$ such that 
for any $\alpha,\beta\in\partial X$ with $\alpha\neq\beta$, 
there exist a sequence $\{g_i\}\subset G$ and a point $y_0\in X$ 
such that $\{d(y_0,g_ix_0)\}_i\rightarrow\infty$ as $i\rightarrow\infty$ and 
for any $i\in\N$, 
\begin{align*}
&\Image \xi_{g_ix_0,\alpha}\cap B(y_0,M)\neq\emptyset \ \text{and} \ \\
&\Image \xi_{g_ix_0,\beta}\cap B(y_0,M)\neq\emptyset,
\end{align*}
where $\xi_{g_ix_0,\alpha}$ and $\xi_{g_ix_0,\beta}$ are 
the geodesic rays with 
$\xi_{g_ix_0,\alpha}(0)=\xi_{g_ix_0,\beta}(0)=g_ix_0$, 
$\xi_{g_ix_0,\alpha}(\infty)=\alpha$ and 
$\xi_{g_ix_0,\beta}(\infty)=\beta$.
Then the boundary $\partial X$ is a scrambled set.
\end{Theorem}

\begin{proof}
By Theorem~\ref{Thm3-1}, 
$$\limsup\{d_{\partial X}(g\alpha,g\beta)\,|\,g\in G\}>0$$ 
for any $\alpha,\beta\in\partial X$ with $\alpha\neq\beta$.
Hence it is sufficient to show that 
$$\liminf\{d_{\partial X}(g\alpha,g\beta)\,|\,g\in G\}=0$$ 
for any $\alpha,\beta\in\partial X$ with $\alpha\neq\beta$.

Let $\alpha,\beta\in\partial X$ with $\alpha\neq\beta$.
By the assumption, 
there exist a sequence $\{g_i\}\subset G$ and a point $y_0\in X$ 
such that $\{d(y_0,g_ix_0)\}_i\rightarrow\infty$ as $i\rightarrow\infty$ and 
\begin{align*}
&\Image \xi_{g_ix_0,\alpha}\cap B(y_0,M)\neq\emptyset \ \text{and} \ \\
&\Image \xi_{g_ix_0,\beta}\cap B(y_0,M)\neq\emptyset
\end{align*}
for any $i\in\N$.
Here we show that 
$$\{d_{\partial X}(g_i^{-1}\alpha,g_i^{-1}\beta)\}_i\rightarrow 0 
\ \text{as}\ i\rightarrow \infty.$$

Let $\epsilon>0$ be a small number and 
let $j_0\in\N$ such that 
$$\frac{1}{2^{j_0-1}}<\epsilon.$$
Since $\{d(y_0,g_ix_0)\}_i\rightarrow\infty$ as $i\rightarrow\infty$, 
there exists $i_0\in\N$ such that for any $i\ge i_0$, 
$$d(y_0,g_ix_0)>\frac{2Mj_0(j_0+1)}{\epsilon}+3M.$$

Let $i\ge i_0$ and let $R_i=d(y_0,g_ix_0)$.
We suppose that 
$\xi_{x_0,g_i^{-1}\alpha}$ and $\xi_{x_0,g_i^{-1}\beta}$ are 
the geodesic rays with 
$\xi_{x_0,g_i^{-1}\alpha}(0)=\xi_{x_0,g_i^{-1}\beta}(0)=x_0$, 
$\xi_{x_0,g_i^{-1}\alpha}(\infty)=g_i^{-1}\alpha$ and 
$\xi_{x_0,g_i^{-1}\beta}(\infty)=g_i^{-1}\beta$.
Then 
$g_i\xi_{x_0,g_i^{-1}\alpha}=\xi_{g_ix_0,\alpha}$ and 
$g_i\xi_{x_0,g_i^{-1}\beta}=\xi_{g_ix_0,\beta}$, 
since $g_i$ is an isometry of $X$.
Here 
\begin{align*}
d_{\partial X}(g_i^{-1}\alpha,g_i^{-1}\beta)
&=\sum_{j=1}^{\infty}\min\{
d(\xi_{x_0,g_i^{-1}\alpha}(j),\xi_{x_0,g_i^{-1}\beta}(j)),\ 
\frac{1}{2^j}\} \\
&\le\sum_{j=1}^{j_0}
d(\xi_{x_0,g_i^{-1}\alpha}(j),\xi_{x_0,g_i^{-1}\beta}(j))
+\sum_{j=j_0+1}^{\infty}\frac{1}{2^j} \\
&=\sum_{j=1}^{j_0}
d(\xi_{g_ix_0,\alpha}(j),\xi_{g_ix_0,\beta}(j))+\frac{1}{2^{j_0}}.
\end{align*}

Now we show that 
$$d(\xi_{g_ix_0,\alpha}(j),\xi_{g_ix_0,\beta}(j))<\frac{2Mj}{R_i-3M}$$
for any $j\in\N$.
By the assumption, 
\begin{align*}
&\Image \xi_{g_ix_0,\alpha}\cap B(y_0,M)\neq\emptyset \ \text{and} \\
&\Image \xi_{g_ix_0,\beta}\cap B(y_0,M)\neq\emptyset.
\end{align*}
Then $\xi_{g_ix_0,\alpha}(t_i)\in B(y_0,M)$ for some $t_i\ge 0$.
Hence $d(y_0,\xi_{g_ix_0,\alpha}(t_i))\le M$.
Here 
\begin{align*}
t_i&=d(g_ix_0,\xi_{g_ix_0,\alpha}(t_i)) \\
&\ge d(g_ix_0,y_0)-d(y_0,\xi_{g_ix_0,\alpha}(t_i)) \\
&\ge d(g_ix_0,y_0)-M \\
&=R_i-M 
\end{align*}
and 
\begin{align*}
d(\xi_{g_ix_0,\alpha}(t_i),\Image \xi_{g_ix_0,\beta})&\le 
d(\xi_{g_ix_0,\alpha}(t_i),y_0)+d(y_0,\Image \xi_{g_ix_0,\beta})\\
&\le 2M.
\end{align*}
By Lemma~\ref{lem1}, 
$$d(\xi_{g_ix_0,\alpha}(t_i-2M),\xi_{g_ix_0,\beta}(t_i-2M))\le 2M.$$
Since $t_i\ge R_i-M$, we have that 
$$d(\xi_{g_ix_0,\alpha}(R_i-3M),\xi_{g_ix_0,\beta}(R_i-3M))\le 2M.$$
Thus by Lemma~\ref{lem0}, 
\begin{align*}
d(\xi_{g_ix_0,\alpha}(j),\xi_{g_ix_0,\beta}(j))
&\le\frac{j}{R_i-3M}
d(\xi_{g_ix_0,\alpha}(R-3M),\xi_{g_ix_0,\beta}(R-3M))\\
&\le \frac{2Mj}{R_i-3M}.
\end{align*}

Hence 
\begin{align*}
d_{\partial X}(g_i^{-1}\alpha,g_i^{-1}\beta)
&\le
\sum_{j=1}^{j_0}
d(\xi_{g_ix_0,\alpha}(j),\xi_{g_ix_0,\beta}(j))+\frac{1}{2^{j_0}} \\
&\le \sum_{j=1}^{j_0}\frac{2Mj}{R_i-3M}+\frac{1}{2^{j_0}} \\
&=\frac{2M}{R_i-3M}\sum_{j=1}^{j_0}j+\frac{1}{2^{j_0}} \\
&=\frac{Mj_0(j_0+1)}{R_i-3M}+\frac{1}{2^{j_0}} \\
&<\frac{\epsilon}{2}+\frac{\epsilon}{2}=\epsilon,
\end{align*}
because 
\begin{align*}
&R_i=d(y_0,g_ix_0)>\frac{2Mj_0(j_0+1)}{\epsilon}+3M \ \text{and}\\
&\displaystyle \frac{1}{2^{j_0-1}}<\epsilon.
\end{align*}

Here $\epsilon>0$ is an arbitrary small number.
Thus 
$$\{d_{\partial X}(g_i^{-1}\alpha,g_i^{-1}\beta)\}_i\rightarrow 0 \ 
\text{as}\ i\rightarrow \infty.$$
This implies that 
$$\liminf\{d_{\partial X}(g\alpha,g\beta)\,|\,g\in G\}=0.$$
Hence the boundary $\partial X$ is a scrambled set.
\end{proof}

\section{Hyperbolic case}

In this section, 
we study boundaries of hyperbolic CAT(0) groups.

We first introduce a definition of hyperbolic spaces.

A geodesic space $X$ is called a {\it hyperbolic space}, 
if there exists a number $\delta>0$ such that 
every geodesic triangle in $X$ is ``$\delta$-thin''.

Here ``$\delta$-thin'' is defined as follows: 
Let $x,y,z \in X$ and 
let $\triangle:=\triangle xyz$ be a geodesic triangle in $X$.
There exist unique non-negative numbers $a,b,c$ such that 
$$d(x,y)=a+b,\ d(y,z)=b+c,\ d(z,x)=c+a.$$
Then we can consider the metric tree $T_\triangle$ that 
has three vertexes of valence one, one vertex of valence three, and 
edges of length $a$, $b$ and $c$.
Let $o$ be the vertex of valence three in $T_\triangle$ and 
let $v_x, v_y, v_z$ be the vertexes of $T_\triangle$ such that 
$d(o,v_x)=a$, $d(o,v_y)=b$ and $d(o,v_z)=c$.
Then the map $\{x,y,z\}\rightarrow\{v_x,v_y,v_z\}$ extends 
uniquely to a map $f:\triangle \rightarrow T_\triangle$ 
whose restriction to each side of $\triangle$ is an isometry.
For some $\delta \ge 0$, 
the geodesic triangle $\triangle$ is said to be {\it $\delta$-thin}, 
if $d(p,q)\le \delta$ for each points $p,q \in \triangle$ with $f(p)=f(q)$.

It is known that 
a geodesic space $X$ is hyperbolic if and only if 
there exists a number $\delta>0$ such that 
every geodesic triangle in $X$ is ``$\delta$-slim''.
Here a geodesic triangle is said to be {\it $\delta$-slim}, 
if each of its sides is contained 
in the $\delta$-neighbourhood of the union of the other two sides.

For a proper hyperbolic space $X$, 
we can define the {\it boundary} $\partial X$ of $X$, 
and if the space $X$ is hyperbolic and CAT(0), then 
these ``boundaries'' coincide.

Details and basic properties of hyperbolic spaces 
and their boundaries are found in 
\cite{BH}, \cite{CP}, \cite{GH} and \cite{Gr}.

A group $G$ is called a {\it hyperbolic group} 
if the group $G$ acts geometrically on some hyperbolic space $X$.
A hyperbolic group $G$ determines the boundary $\partial X$ 
of a hyperbolic space $X$ on which $G$ acts geometrically.
The boundary $\partial X$ is called the {\it boundary of} $G$ 
and denoted by $\partial G$.

It is known when is a CAT(0) space hyperbolic.

\begin{Theorem}[{\cite[p.400, Theorem~III.H.1.5]{BH}}]\label{Thm:Flat}
A proper cocompact CAT(0) space $X$ is hyperbolic if and only if 
it does not contain a subspace isometric to the flat plane $\R^2$.
\end{Theorem}

We show that the boundary of 
a non-elementary hyperbolic CAT(0) group 
is always a scrambled set.

\begin{Theorem}\label{Thm:Hy}
Suppose that a group $G$ acts 
geometrically on a CAT(0) space $X$ and $|\partial X|>2$.
If $X$ is hyperbolic, 
then the boundary $\partial X$ is a scrambled set.
Moreover there exists a constant $c>0$ such that 
$$\limsup\{d_{\partial X}(g\alpha,g\beta)\,|\,g\in G\}>c$$ 
for any $\alpha,\beta\in\partial X$ with $\alpha\neq\beta$.
\end{Theorem}

\begin{proof}
Suppose that $X$ is hyperbolic.
There exists a constant $\delta>0$ such that 
every geodesic triangle in $X$ is $\delta$-thin.
Since the action of $G$ on $X$ is cocompact and $X$ is proper, 
$GB(x_0,N)=X$ for some $N>0$.

Let $\alpha,\beta\in\partial X$ with $\alpha\neq\beta$.
Since $|\partial X|>2$, 
there exists $\gamma\in\partial X\setminus\{\alpha,\beta\}$.
Let $x_i=\xi_{x_0,\gamma}(i)$ for each $i\in\N$, 
where $\xi_{x_0,\gamma}$ is the geodesic ray 
with $\xi_{x_0,\gamma}(0)=x_0$ and $\xi_{x_0,\gamma}(\infty)=\gamma$.
We consider the two triangles of $x_0,x_i,\alpha$ and $x_0,x_i,\beta$.
Since $X$ is hyperbolic, 
there exists an enough large number $R>0$ such that for any $i\in\N$ with $i>R$, 
\begin{align*}
&d(\xi_{x_0,\gamma}(R),\Image \xi_{x_i,\alpha})\le\delta\ \text{and}\ \\
&d(\xi_{x_0,\gamma}(R),\Image \xi_{x_i,\beta})\le\delta.
\end{align*}
For each $i\in\N$ with $i>R$, 
there exists $g_i\in G$ such that $d(x_i,g_ix_0)\le N$ 
by the definition of the number $N>0$.
Then 
$d(\xi_{x_0,\gamma}(R),\Image \xi_{x_i,\alpha})\le\delta$ and 
$d(\xi_{x_i,\alpha}(t),\xi_{g_ix_0,\alpha}(t))\le N$ for any $t\ge 0$ 
by Lemma~\ref{lem00}~(2).
Hence 
$$d(\xi_{x_0,\gamma}(R),\Image \xi_{g_ix_0,\alpha})\le N+\delta.$$
Also we have that 
$$d(\xi_{x_0,\gamma}(R),\Image \xi_{g_ix_0,\beta})\le N+\delta.$$
Thus the constant $M=N+\delta$, 
the point $y_0=\xi_{x_0,\gamma}(R)$ and 
the sequence $\{g_i\,|\,i\in\N, i>R\}$ 
satisfy the condition of Theorem~\ref{Thm:key}.
Hence the boundary $\partial X$ is a scrambled set.

Moreover, 
for each $\alpha,\beta\in\partial X$ with $\alpha\neq\beta$, 
$\angle(\alpha,\beta)=\pi$ 
because there exists a geodesic line $\sigma:\R\rightarrow X$ 
such that $\sigma(\infty)=\alpha$ and $\sigma(-\infty)=\beta$ 
(cf.\ \cite[p.428, Lemma~III.H.3.2]{BH} and \cite{CP}).
Hence by the proof of Theorem~\ref{Thm3-1}, 
for the constants $r=2\sin(\pi/4)$ and 
$\displaystyle t_0=[\frac{2N+1}{r}]+1$, 
there exists a sequence $\{g_i\}\subset G$ such that 
$$ d_{\partial X}(g_i^{-1}\alpha,g_i^{-1}\beta) 
\ge\frac{1}{2^{t_0}}>\frac{1}{2^{t_0+1}}. $$
Here $c:=\displaystyle \frac{1}{2^{t_0+1}}$ is a constant 
which does not depend on $\alpha$ and $\beta$.
Thus 
$$\limsup\{d_{\partial X}(g\alpha,g\beta)\,|\,g\in G\}>c$$ 
for any $\alpha,\beta\in\partial X$ with $\alpha\neq\beta$.
\end{proof}

It is known that 
the boundary $\partial G$ of a hyperbolic group $G$ is minimal.
Hence the boundary $\partial G$ 
of every non-elementary hyperbolic CAT(0) group $G$ is 
a scrambled set and minimal.

Also we obtain a theorem.

\begin{Theorem}\label{Thm:Hy2}
Suppose that a group $G$ acts 
geometrically on a CAT(0) space $X$ and $|\partial X|>2$.
The space $X$ is hyperbolic if and only if 
there exists a constant $c>0$ such that 
$$\limsup\{d_{\partial X}(g\alpha,g\beta)\,|\,g\in G\}>c$$ 
for any $\alpha,\beta\in\partial X$ with $\alpha\neq\beta$.
\end{Theorem}

\begin{proof}
By Theorem~\ref{Thm:Hy}, 
it is sufficient to show that 
if $X$ is not hyperbolic then 
there does not exist a constant $c>0$ such that 
$$\limsup\{d_{\partial X}(g\alpha,g\beta)\,|\,g\in G\}>c$$ 
for any $\alpha,\beta\in\partial X$ with $\alpha\neq\beta$.

Suppose that $X$ is not hyperbolic.
Then $X$ contsins a subspace $Z$ which is isometric to the flat plane $\R^2$ 
by Theorem~\ref{Thm:Flat}.
The boundary of $Z$ is the circle, and 
for each $\theta\in[0,\pi]$ 
there exist $\alpha,\beta\in\partial Z$ 
with $\angle(\alpha,\beta)=\theta$.

Let $\epsilon>0$ be a small number.
There exists $t_0\in\N$ such that 
$\displaystyle \frac{1}{2^{t_0}}<\frac{\epsilon}{2}$.
Then there exist $\alpha,\beta\in\partial Z\subset\partial X$ with 
$\alpha\neq\beta$ such that 
$$\sin(\angle(\alpha,\beta)/2)\le \frac{\epsilon}{4t_0^2}.$$
For each $y_0\in X$, 
Lemmas~\ref{lem0} and \ref{lem:angle}~(2) imply that 
\begin{align*}
\frac{1}{t_0}d(\xi_{y_0,\alpha}(t_0),\xi_{y_0,\beta}(t_0))
&\le \lim_{t\rightarrow\infty}\frac{1}{t}d(\xi_{y_0,\alpha}(t),\xi_{y_0,\beta}(t)) \\
&=2\sin(\angle(\alpha,\beta)/2) \\
&\le 2(\frac{\epsilon}{4t_0^2})=\frac{\epsilon}{2t_0^2},
\end{align*}
where $\xi_{y_0,\alpha}$ and $\xi_{y_0,\beta}$ are the geodesic rays 
with $\xi_{y_0,\alpha}(0)=\xi_{y_0,\beta}(0)=y_0$, 
$\xi_{y_0,\alpha}(\infty)=\alpha$ and $\xi_{y_0,\beta}(\infty)=\beta$.
Hence 
$$d(\xi_{y_0,\alpha}(t_0),\xi_{y_0,\beta}(t_0))
\le \frac{\epsilon}{2t_0}$$
for any $y_0\in X$.
Thus for each $j\in\{1,2,\dots,t_0\}$, 
$$d(\xi_{y_0,\alpha}(j),\xi_{y_0,\beta}(j))
\le d(\xi_{y_0,\alpha}(t_0),\xi_{y_0,\beta}(t_0))
\le \frac{\epsilon}{2t_0}$$
by Lemma~\ref{lem0}.

By the above argument, for each $g\in G$, 
\begin{align*}
d_{\partial X}(g\alpha,g\beta)
&=\sum_{j=1}^{\infty}\min\{
d(\xi_{x_0,g\alpha}(j),\xi_{x_0,g\beta}(j)),\ \frac{1}{2^j}\} \\
&=\sum_{j=1}^{\infty}\min\{
d(\xi_{g^{-1}x_0,\alpha}(j),\xi_{g^{-1}x_0,\beta}(j)),\ \frac{1}{2^j}\} \\
&\le
\sum_{j=1}^{t_0}
d(\xi_{g^{-1}x_0,\alpha}(j),\xi_{g^{-1}x_0,\beta}(j))
+\sum_{j=t_0+1}^{\infty}\frac{1}{2^j} \\
&\le \sum_{j=1}^{t_0}\frac{\epsilon}{2t_0}+\frac{1}{2^{t_0}} \\
&\le\frac{\epsilon}{2}+\frac{\epsilon}{2}=\epsilon.
\end{align*}

Thus 
for any $\epsilon>0$ 
there exists $\alpha,\beta\in\partial X$ with $\alpha\neq\beta$ 
such that 
$$\limsup\{d_{\partial X}(g\alpha,g\beta)\,|\,g\in G\}\le\epsilon.$$
This means that 
there does not exist a constant $c>0$ such that 
$$\limsup\{d_{\partial X}(g\alpha,g\beta)\,|\,g\in G\}>c$$ 
for any $\alpha,\beta\in\partial X$ with $\alpha\neq\beta$.
\end{proof}

An action of a group $G$ on a metric space $Y$ by homeomorphisms 
is said to be {\it expansive}, 
if there exists a constant $c>0$ such that 
for each pair $y,y'\in Y$ with $y\neq y'$, 
there is $g\in G$ such that $d(gy,gy')>c$.

We obtain a corollary from the proofs of Theorems~\ref{Thm:Hy} and \ref{Thm:Hy2}.

\begin{Corollary}\label{Cor:Hy3}
Suppose that a group $G$ acts 
geometrically on a CAT(0) space $X$ and $|\partial X|>2$.
The action of $G$ on $\partial X$ is expansive if and only if 
the space $X$ is hyperbolic.
\end{Corollary}

\section{CAT(0) groups whose boundaries are scrambled sets}

In this section, 
we investigate when are boundaries of CAT(0) groups scrambled sets.

We first give a sufficient condition of 
CAT(0) groups whose boundaries are scrambled sets.

Suppose that a group $G$ acts geometrically on a CAT(0) space $X$.
For an element $g\in G$, 
we define $Z_g$ as the centralizer of $g$ and 
define $F_{g}$ as the fixed-point set of $g$ in $X$, 
that is, 
\begin{align*}
Z_g&=\{h\in G|\, gh=hg\} \ \text{and}\\
F_{g}&=\{x\in X|\, gx=x\}.  \\
\end{align*}

The following lemmas are known.

\begin{Lemma}[{\cite[Theorem~2.1]{Ho3} and \cite{Ru}}]\label{lem:fix}
Suppose that a group $G$ acts geometrically on a CAT(0) space $X$.
For $g\in G$, if $Z_g$ is finite then $F_g$ is bounded.
\end{Lemma}

\begin{Lemma}[{\cite[Proposition~II.6.2(2)]{BH}}]\label{lem:fix2}
Let $X$ be a CAT(0) space and 
let $g$ and $h$ be isometries of $X$.
Then $gF_h=F_{ghg^{-1}}$.
\end{Lemma}

Here we show a theorem.

\begin{Theorem}\label{Thm:CAT(0)}
Suppose that a group $G$ acts geometrically on a CAT(0) space $X$.
If there exists an element $g_0\in G$ such that 
\begin{enumerate}
\item[(1)] $Z_{g_0}$ is finite, 
\item[(2)] $X\setminus F_{g_0}$ is not connected, and 
\item[(3)] each component of $X\setminus F_{g_0}$ is convex and not $g_0$-invariant, 
\end{enumerate}
then the boundary $\partial X$ is a scrambled set.
\end{Theorem}

\begin{proof}
Suppose that there exists an element $g_0\in G$ 
which satisfies the conditions (1), (2) and (3).
Using Theorem~\ref{Thm:key}, 
we show that the boundary $\partial X$ is a scrambled set.

Let $\alpha,\beta\in\partial X$ with $\alpha\neq\beta$.

By (1) and Lemma~\ref{lem:fix}, 
$F_{g_0}$ is bounded.
Hence for some (and any) $y_0\in F_{g_0}$, 
$$L(\bigcup\{gF_{g_0}\,|\,g\in G\})=L(Gy_0)=\partial X. $$
Here for a subset $A\subset X$, 
the limit set $L(A)$ of $A$ is defined as 
$$L(A)=\overline{A}\cap \partial X, $$
where $\overline{A}$ is the closure of $A$ in $X\cup \partial X$.

Now $|\partial X|>2$ (hence $\partial X$ is uncoundatble).
Then there exists $h\in G$ such that 
$$(\Image\xi_{x_0,\alpha}\cup\Image\xi_{x_0,\beta})\cap hF_{g_0}=\emptyset.$$
Indeed if 
$$(\Image\xi_{x_0,\alpha}\cup\Image\xi_{x_0,\beta})\cap gF_{g_0}\neq\emptyset$$ 
for any $g\in G$, then 
$$\{\alpha,\beta\}=L(\Image\xi_{x_0,\alpha}\cup\Image\xi_{x_0,\beta})=\partial X,$$
which contradicts to $|\partial X|>2$.

Let $h\in G$ such that 
$$(\Image\xi_{x_0,\alpha}\cup\Image\xi_{x_0,\beta})\cap hF_{g_0}=\emptyset.$$
By Lemma~\ref{lem:fix2}, $hF_{g_0}=F_{hg_0h^{-1}}$.
We consider the point $hg_0h^{-1}x_0$.
By (2) and (3), 
$X\setminus hF_{g_0}$ is not connected, and 
two points $x_0$ and $hg_0h^{-1}x_0$ are 
in distinct components of $X\setminus hF_{g_0}$.
Let $A_0$ and $A_1$ be the components of $X\setminus hF_{g_0}$ 
such that $x_0\in A_0$ and $hg_0h^{-1}x_0\in A_1$.
Then $\Image\xi_{x_0,\alpha}\cup\Image\xi_{x_0,\beta}\subset A_0$.
Since the components $A_0$ and $A_1$ are unbounded, 
we can take a sequence $\{h_i\}\subset G$ 
such that $\{h_ix_0\}_i\subset A_1$ and 
$\{d(h_ix_0,hF_{g_0})\}_i\rightarrow\infty$ as $i\rightarrow\infty$.
Since $h_ix_0\in A_1$ and $\alpha,\beta\in\partial A_0$, 
\begin{align*}
&\Image \xi_{h_ix_0,\alpha}\cap hF_{g_0}\neq\emptyset \ \text{and} \\
&\Image \xi_{h_ix_0,\beta}\cap hF_{g_0}\neq\emptyset 
\end{align*}
for each $i$.
For the diameter $M=\diam(F_{g_0})=\diam(hF_{g_0})$ 
and some point $y_0\in hF_{g_0}$, 
we have that 
\begin{align*}
&\Image \xi_{h_ix_0,\alpha}\cap B(y_0,M)\neq\emptyset \ \text{and} \\
&\Image \xi_{h_ix_0,\beta}\cap B(y_0,M)\neq\emptyset. 
\end{align*}
Here $M>0$ is a constant which does not depend on $\alpha$ and $\beta$.

Thus the condition of Theorem~\ref{Thm:key} holds, 
and the boundary $\partial X$ is a scrambled set.
\end{proof}

It is known that 
if the condition of Theorem~\ref{Thm:CAT(0)} holds 
then the boundary $\partial X$ is also minimal (\cite{Ho3}).

In this paper, 
we define a {\it reflection} of a geodesic space as follows:
An isometry $r$ of a geodesic space $X$ is called a 
{\it reflection} of $X$, if 
\begin{enumerate}
\item[(1)] $r^2$ is the identity of $X$, 
\item[(2)] $X\setminus F_r$ has exactly 
two convex connected components $X^+_r$ and $X^-_r$ and 
\item[(3)] $rX^+_r=X^-_r$, 
\end{enumerate}
where $F_r$ is the fixed-points set of $r$.
We note that ``reflections'' in this paper need not 
satisfy the condition (4) $\Int F_r=\emptyset$ in \cite{Ho2-2}.

We obtain a corollary from Theorem~\ref{Thm:CAT(0)}.

\begin{Corollary}\label{Cor:CAT(0)}
Suppose that a group $G$ acts geometrically on a CAT(0) space $X$ 
and $|\partial X|>2$.
If there exists a reflection $r\in G$ of $X$ such that $Z_r$ is finite, 
then the boundary $\partial X$ is a scrambled set.
\end{Corollary}

Next, we give a sufficient condition of CAT(0) groups 
whose boundaries are {\it not} scrambled sets.

A subset $A$ of a metric space $Y$ is said to be {\it quasi-dense}, 
if there exists a constant $K>0$ such that 
for each $y\in Y$ there is $a\in A$ such that $d(y,a)\le K$, 
that is, 
$$B(A,K)=\{y\in Y\,|\,d(y,a)\le K\ \text{for some}\ a\in A\}=Y.$$

\begin{Theorem}\label{Thm:non}
Suppose that a group $G$ acts 
geometrically on a CAT(0) space $X$ and $|\partial X|>2$.
If $X$ contains a quasi-dense subset $X_1\times X_2$ 
such that $X_1$ and $X_2$ are unbounded, 
then the boundary $\partial X$ is not a scrambled set.
\end{Theorem}

\begin{proof}
Suppose that $X$ contains a quasi-dense subset $X_1\times X_2$ 
such that $X_1$ and $X_2$ are unbounded.
Then there exists a constant $K>0$ 
such that 
$$B(X_1\times X_2,K)=X.$$
Since $X_1$ and $X_2$ are unbounded, 
there exist $\alpha\in\partial X_1$ and $\beta\in\partial X_2$.
We note that 
$$\partial X=\partial(X_1\times X_2)=\partial X_1*\partial X_2,$$
where $\partial X_1*\partial X_2$ is the spherical join.
Hence $\angle(\alpha,\beta)=\pi/2$ 
and $\angle_{z_0}(\alpha,\beta)=\pi/2$ for each $z_0\in X_1\times X_2$.

To show that the boundary $\partial X$ is not a scrambled set, 
we prove that 
$$\liminf\{d_{\partial X}(g\alpha,g\beta)\,|\,g\in G\}>0,$$
that is, 
we show that there exists a constant $c>0$ such that 
$d_{\partial X}(g\alpha,g\beta)\ge c$ for any $g\in G$.

Let $\displaystyle t_0=[\frac{2K+1}{\sqrt{2}}]+1$ and 
$\displaystyle c=\frac{1}{2^{t_0}}$.
Then we show that 
$d_{\partial X}(g\alpha,g\beta)\ge c$ for any $g\in G$.

Let $g\in G$.
Since $B(X_1\times X_2,K)=X$, 
there exists $z_0\in X_1\times X_2$ such that $d(g^{-1}x_0,z_0)\le K$.
We consider the geodesic rays 
$\xi_{z_0,\alpha}$, $\xi_{z_0,\beta}$, 
$\xi_{g^{-1}x_0,\alpha}$ and $\xi_{g^{-1}x_0,\beta}$.
By Lemma~\ref{lem00}~(2), 
\begin{align*}
&d(\xi_{g^{-1}x_0,\alpha}(t),\xi_{z_0,\alpha}(t))\le K \ \text{and} \\
&d(\xi_{g^{-1}x_0,\beta}(t),\xi_{z_0,\beta}(t))\le K.
\end{align*}
Since $\angle_{z_0}(\alpha,\beta)=\pi/2$ and 
the convex hull of $\Image\xi_{z_0,\alpha}\cup\Image\xi_{z_0,\beta}$ is flat, 
we have that $d(\xi_{z_0,\alpha}(t),\xi_{z_0,\beta}(t))=\sqrt{2}t$ 
for any $t\ge 0$.
Hence 
\begin{align*}
d(\xi_{g^{-1}x_0,\alpha}(t),\xi_{g^{-1}x_0,\beta}(t))
&\ge d(\xi_{z_0,\alpha}(t),\xi_{z_0,\beta}(t))-2K \\
&= \sqrt{2}t-2K.
\end{align*}
Here $\sqrt{2}t-2K\ge 1$ if and only if 
$\displaystyle t\ge \frac{2K+1}{\sqrt{2}}$.
Hence for $\displaystyle t_0=[\frac{2K+1}{\sqrt{2}}]+1$, 
we obtain that 
$$d(\xi_{g^{-1}x_0,\alpha}(t_0),\xi_{g^{-1}x_0,\beta}(t_0))
\ge\sqrt{2}t_0-2K\ge 1.$$
Then 
\begin{align*}
d_{\partial X}(g\alpha,g\beta)
&=\sum_{j=1}^{\infty}\min\{
d(\xi_{x_0,g\alpha}(j),\xi_{x_0,g\beta}(j)),\ 
\frac{1}{2^j}\} \\
&=\sum_{j=1}^{\infty}\min\{
d(\xi_{g^{-1}x_0,\alpha}(j),\xi_{g^{-1}x_0,\beta}(j)),\ 
\frac{1}{2^j}\} \\
&\ge \frac{1}{2^{t_0}}=c.
\end{align*}
Here $g\in G$ is arbitrary.
Thus 
$$\liminf\{d_{\partial X}(g\alpha,g\beta)\,|\,g\in G\}\ge c>0, $$
and the boundary $\partial X$ is not a scrambled set.
\end{proof}

By a splitting theorem on CAT(0) spaces 
(\cite{Ho4} and \cite{Mo}, cf.\ \cite{Ho5}), 
we obtain that 
if a CAT(0) group $G$ contains a subgroup $G_1\times G_2$ of finite index 
such that $G_1$ and $G_2$ are infinite, 
then a CAT(0) space $X$ on which $G$ acts geometrically 
contains a quasi-dense subspace which splits as a product $X_1\times X_2$, 
where $X_1$ and $X_2$ are unbounded.
This implies the following corollary.

\begin{Corollary}\label{Cor:non}
Suppose that a group $G$ acts geometrically on a CAT(0) space $X$ 
and $|\partial X|>2$.
If $G$ contains a subgroup $G_1\times G_2$ of finite index 
such that $G_1$ and $G_2$ are infinite, 
then the boundary $\partial X$ is not a scrambled set.
\end{Corollary}

\section{Coxeter groups and Davis complexes}

In this section, 
we introduce definitions and some properties of 
Coxeter groups and Davis complexes.

A {\it Coxeter group} is a group $W$ having a presentation
$$\langle \,S \, | \, (st)^{m(s,t)}=1 \ \text{for}\ s,t \in S \,
\rangle,$$ 
where $S$ is a finite set and 
$m:S \times S \rightarrow \N \cup \{\infty\}$ is a function 
satisfying the following conditions:
\begin{enumerate}
\item[(1)] $m(s,t)=m(t,s)$ for any $s,t \in S$,
\item[(2)] $m(s,s)=1$ for any $s \in S$, and
\item[(3)] $m(s,t) \ge 2$ for any $s,t \in S$ with $s\neq t$.
\end{enumerate}
The pair $(W,S)$ is called a {\it Coxeter system}.
If, in addition, 
\begin{enumerate}
\item[(4)] $m(s,t) = 2$ or $\infty$ for any $s,t \in S$ with $s\neq t$,
\end{enumerate}
then $(W,S)$ is said to be {\it right-angled}.
We note that for $s,t\in S$, 
$m(s,t)=2$ if and only if $st=ts$.

Let $(W,S)$ be a Coxeter system.
For $w\in W$, 
we denote by $\ell(w)$ the word length of $w$ with respect to $S$.
For $w\in W$, 
a representation $w=s_1\cdots s_l$ ($s_i \in S$) is said to be 
{\it reduced}, if $\ell(w)=l$.
The Coxeter group $W$ has the {\it word metric} $d_{\ell}$ 
defined by $d_{\ell}(w,w')=\ell(w^{-1}w')$ for each $w,w'\in W$.
Also for a subset $T \subset S$, 
$W_T$ is defined as the subgroup of $W$ generated by $T$, 
and called a {\it parabolic subgroup}.
It is known that the pair $(W_T,T)$ is also a Coxeter system (\cite{Bo}).
If $T$ is the empty set, then $W_T$ is the trivial group.
A subset $T\subset S$ is called a {\it spherical subset of $S$}, 
if the parabolic subgroup $W_T$ is finite.

Let $(W,S)$ be a Coxeter system.
For each $w \in W$, we define a subset $S(w)$ of $S$ as 
$$ S(w)=\{s \in S \,|\, \ell(ws) < \ell(w)\}.$$
Also for a subset $T$ of $S$, we define a subset $W^T$ of $W$ as 
$$ W^T=\{w \in W \,|\, S(w)=T\}.$$

The following lemma is known.

\begin{Lemma}[\cite{Bo}, \cite{D1}, \cite{D3}]\label{lem:Coxeter-basic}
Let $(W,S)$ be a Coxeter system.
For each $w\in W$, 
$S(w)$ is a spherical subset of $S$, 
i.e., $W_{S(w)}$ is finite.
\end{Lemma}

Every Coxeter system $(W,S)$ determines 
a {\it Davis complex} $\Sigma(W,S)$ 
which is a CAT(0) geodesic space (\cite{D1}, \cite{D2}, \cite{D3}, \cite{M}).
Here the vertex set of $\Sigma(W,S)$ is $W$ and 
the $1$-skeleton of $\Sigma(W,S)$ is 
the Cayley graph of $W$ with respect to $S$.
Also $\Sigma(W,S)$ is contractible.
The natural action of $W$ on $\Sigma(W,S)$ 
is proper, cocompact and by isometries, i.e., 
the Coxeter group $W$ acts geometrically on the Davis complex $\Sigma(W,S)$ and 
$W$ is a CAT(0) group.
If $W$ is infinite, then $\Sigma(W,S)$ is noncompact and 
we can consider the boundary $\partial\Sigma(W,S)$ of 
the CAT(0) space $\Sigma(W,S)$.
This boundary 
$\partial\Sigma(W,S)$ is called the {\it boundary of} $(W,S)$.

Let $(W,S)$ be a Coxeter system.
The set $R_S=\{wsw^{-1}\,|\, w\in W, s\in S\}$ is called 
the set of {\it reflections} of $(W,S)$.
In fact, each $r\in R_S$ is a ``reflection'' of the Davis complex $\Sigma(W,S)$.
We define $K(W,S)$ as the closure of the component $C$ 
of $\Sigma(W,S)\setminus \bigcup_{r\in R_S}F_r$ with $1\in C$, 
where $F_r$ is the fixed-point set of $r$ in $\Sigma(W,S)$.
It is known that the subset $K(W,S)$ is compact and 
a fundamental domain of the action of $W$ on $\Sigma(W,S)$, 
that is, $WK(W,S)=\Sigma(W,S)$.

The following lemmas are known.
These lemmas give a relation between geodesic paths in Davis complexes 
and reduced words in Coxeter systems.

\begin{Lemma}[{\cite[Lemma~4.2]{Ho0}}]\label{Lem:word1}
Let $(W,S)$ be a Coxeter system and 
let $N$ be the diameter of $K(W,S)$ in $\Sigma(W,S)$.
Then 
for any $w \in W$ with $w\neq 1$, 
there exists a reduced representation $w=s_1\cdots s_l$ such that 
$$d(s_1\cdots s_i,[1,w])\le N$$
for any $i\in\{1,\dots,l\}$, 
where $[1,w]$ is the geodesic from $1$ to $w$ in $\Sigma(W,S)$.
\end{Lemma}

\begin{Lemma}[{\cite[Lemma~2.6]{Ho2}}]\label{Lem:word2}
Let $(W,S)$ be a Coxeter system and 
let $N$ be the diameter of $K(W,S)$ in $\Sigma(W,S)$.
Then for any $\alpha\in\partial\Sigma(W,S)$ 
there exists a sequense $\{s_i\}\subset S$ 
such that $s_1\cdots s_i$ is reduced and 
$$d(s_1\cdots s_i,\Image\xi_{1,\alpha})\le N$$ 
for any $i\in\N$, 
where 
$\xi_{1,\alpha}$ is the geodesic ray in $\Sigma(W,S)$ with 
$\xi_{1,\alpha}(0)=1$ and $\xi_{1,\alpha}(\infty)=\alpha$.
\end{Lemma}

Also there is the following lemma.

\begin{Lemma}[{\cite[Lemma~3.3]{Ho2}}]\label{Lem:word3}
Let $(W,S)$ be a Coxeter system, 
let $N$ be the diameter of $K(W,S)$ in $\Sigma(W,S)$ 
and let $x,y\in W$.
If $o(st)=\infty$ for each $s\in S(x)$ and $t\in S(y^{-1})$, 
then $d(x,[1,xy]) \le N$, 
where $o(st)$ is the order of $st$ in $W$.
\end{Lemma}

\section{Coxeter systems whose boundaries are scrambled sets}

In this section, 
we investigate Coxeter systems whose boundaries are scrambled sets.
We give sufficient conditions of a Coxeter system 
whose boundary is a scrambed set.

\begin{Theorem}\label{Thm7-1}
Let $(W,S)$ be a Coxeter system with $|\partial\Sigma(W,S)|>2$.
Suppose that there exist $s_0,t_0\in S$ and a number $K>0$ such that 
\begin{enumerate}
\item[(1)] $o(s_0t_0)=\infty$ and 
\item[(2)] for each $w,v\in W$, 
there exists $x\in W$ such that 
$\ell(x)\le K$ and $wx,vx\in W^{\{s_0\}}$.
\end{enumerate}
Then the boundary $\partial\Sigma(W,S)$ is a scrambled set.
\end{Theorem}

\begin{proof}
Let $\alpha,\beta\in\partial\Sigma(W,S)$ with $\alpha\neq\beta$.
By Lemma~\ref{Lem:word2}, 
there exist sequences $\{a_i\},\{b_i\}\subset S$ 
such that 
\begin{align*}
&d(a_1\cdots a_i,\Image\xi_{1,\alpha})\le N \ \text{and} \\
&d(b_1\cdots b_i,\Image\xi_{1,\beta})\le N
\end{align*}
for each $i\in\N$, where $N=\diam(K(W,S))$.
Let $w_i=a_1\cdots a_i$ and $v_i=b_1\cdots b_i$ for each $i$.
By (2), for each $i$, 
there exists $x_i\in W$ such that 
$\ell(x_i)\le K$ and $w_i^{-1}x_i,v_i^{-1}x_i\in W^{\{s_0\}}$.
Since $\{x\in W\,|\,\ell(x)\le K\}$ is finite, 
there exist $x\in W$ and a subsequence $\{i_j\,|\,j\in\N\}\subset\N$ 
such that $x_{i_j}=x$ for any $j\in\N$.
We note that 
the sequences $\{x^{-1}w_{i_j}\}$ and $\{x^{-1}v_{i_j}\}$ 
converge to $x^{-1}\alpha$ and $x^{-1}\beta$ respectively, 
$(x^{-1}w_{i_j})^{-1}=w_{i_j}^{-1}x\in W^{\{s_0\}}$ and 
$(x^{-1}v_{i_j})^{-1}=v_{i_j}^{-1}x\in W^{\{s_0\}}$.
By Lemma~\ref{Lem:word3}, 
$$d((s_0t_0)^k,[1,(s_0t_0)^kx^{-1}w_{i_j}])\le N$$
for each $k\in\N$ and $j\in\N$, 
because $o(s_0t_0)=\infty$ by (1) and $(s_0t_0)^k\in W^{\{t_0\}}$.
Hence 
$$d((s_0t_0)^k,\Image\xi_{1,(s_0t_0)^kx^{-1}\alpha})\le N$$
for each $k\in\N$.
Let $g_k=(s_0t_0)^kx^{-1}$ for each $k\in\N$.
Since $g_k^{-1}$ is an isometry of $\Sigma(W,S)$, 
we have that 
$$d(g_k^{-1}(s_0t_0)^k,
\Image\xi_{g_k^{-1},g_k^{-1}(s_0t_0)^kx^{-1}\alpha})\le N.$$
Hence 
$$d(x,\Image\xi_{g_k^{-1},\alpha})\le N$$
for each $k\in\N$.
By the same argument, we also obtain that 
$$d(x,\Image\xi_{g_k^{-1},\beta})\le N$$
for each $k\in\N$.
Here $d(1,x)\le K$ because $\ell(x)\le K$.
Hence 
\begin{align*}
&d(1,\Image\xi_{g_k^{-1},\alpha})\le N+K \ \text{and} \\
&d(1,\Image\xi_{g_k^{-1},\beta})\le N+K 
\end{align*}
for each $k\in\N$.
We note that $\{d(1,g_k^{-1})\}_k\rightarrow\infty$ as $k\rightarrow\infty$ 
and the number $M:=N+K$ is a constant 
which does not depend on $\alpha$ and $\beta$.

Thus the condition of Theorem~\ref{Thm:key} holds, 
and the boundary $\partial\Sigma(W,S)$ is a scrambled set 
\end{proof}

There is a similarity between the condiotions of 
Theorem~\ref{Thm7-1} and \cite[Theorem~3.1]{Ho5} (and \cite[Theorem~4.1]{Ho1}).
Here \cite[Theorem~3.1]{Ho5} and \cite[Theorem~4.1]{Ho1} 
give a sufficient condition of Coxeter systems whose boundaries are minimal.

For a Coxeter system $(W,S)$ and the Davis complex $\Sigma(W,S)$, 
each $s\in S$ is a ``reflection'' of $\Sigma(W,S)$ 
in the sence of the definition in Section~5.
Hence we obtain a corollary from Corollary~\ref{Cor:CAT(0)}.

\begin{Corollary}\label{Cor:Coxeter2}
Let $(W,S)$ be a Coxeter system with $|\partial\Sigma(W,S)|>2$.
If there exists $s\in S$ such that $Z_s$ is finite, 
then the boundary $\partial\Sigma(W,S)$ is a scrambled set.
\end{Corollary}

We give an example of a Coxeter group 
which is not hyperbolic and 
which satisfies the condition in Corollary~\ref{Cor:Coxeter2} 
(hence this example satisfies the conditions in 
Theorem~\ref{Thm:CAT(0)} and Corollary~\ref{Cor:CAT(0)}).

\begin{Example}
We consider the Coxeter system $(W,S)$ defined by 
the diagram in Figure~1.
Since $\Sigma(W_{\{t_1,t_2,t_3\}},\{t_1,t_2,t_3\})$ is 
the flat Euclidean plane, 
the Coxeter group $W$ is not hyperbolic.
Also $Z_s=W_{\{s,t_1,t_2\}}$ is finite.
Hence $(W,S)$ satisfies the condition of Corollary~\ref{Cor:Coxeter2}, 
and the boundary $\partial\Sigma(W,S)$ is a scrambled set.
\begin{figure}[htbp]
\unitlength = 0.9mm
\begin{center}
\begin{picture}(80,30)(-40,-5)
\put(0,0){\line(-3,2){18}}
\put(0,0){\line(3,2){18}}
\put(0,0){\line(0,1){24}}
\put(0,24){\line(-3,-2){18}}
\put(0,24){\line(3,-2){18}}
\put(0,0){\circle*{1.3}}
\put(0,24){\circle*{1.3}}
\put(-18,12){\circle*{1.3}}
\put(18,12){\circle*{1.3}}
{\small
\put(1.5,-3){$t_2$}
\put(18.5,9){$s$}
\put(1,24){$t_1$}
\put(-20,9){$t_3$}
\put(-10.5,2.5){$2$}
\put(-10.5,19){$4$}
\put(9.5,3){$2$}
\put(9.5,19){$2$}
\put(1,11){$4$}
}
\end{picture}
\end{center}
\caption[]{}
\label{fig1}
\end{figure}
\end{Example}

\section{Right-angled Coxeter groups}

The purpose of this section is to prove the following theorem.

\begin{Theorem}\label{Thm8-1}
If $(W,S)$ is an irreducible right-angled Coxeter system 
and $|\partial\Sigma(W,S)|>2$, 
then the boundary $\partial\Sigma(W,S)$ is a scrambled set.
\end{Theorem}

A Coxeter system $(W,S)$ is said to be {\it irreducible} 
if, for any nonempty and proper subset $T$ of $S$, 
$W$ does not decompose as 
the direct product of $W_T$ and $W_{S \setminus T}$.

A Coxeter group $W$ is said to be {\it right-angled}, 
if $(W,S)$ is a right-angled Coxeter system for some $S\subset W$.
It is known that every right-angled Coxeter group determines its Coxeter system 
up to isomorphisms (\cite{R}).
Hence a right-angled Coxeter group $W$ 
determines the boundary $\partial\Sigma(W,S)$.

The following lemmas are known.

\begin{Lemma}[\cite{Bo}, \cite{Hu}]\label{lem:right-angled:basic}
Let $(W,S)$ be a right-angled Coxeter system.
\begin{enumerate}
\item[(1)] $W$ is finite if and only if $st=ts$ for any $s,t\in S$, 
that is, $W$ is isomorphic to $(\Z_2)^{|S|}$.
\item[(2)] $(W,S)$ is irreducible if and only if 
for each $a,b\in S$ with $a\neq b$ 
there exists a sequence $\{a=s_1,s_2,\dots,s_n=b\}\subset S$ 
such that $o(s_is_{i+1})=\infty$ for any $i\in\{1,\dots,n-1\}$.
\end{enumerate}
\end{Lemma}

\begin{Lemma}[{\cite[Lemma~2.7]{Ho5}}]\label{Lem:right-angled:2}
Let $(W,S)$ be a right-angled Coxeter system, 
let $U$ be a spherical subset of $S$, 
let $s_0\in S\setminus U$ and 
let $T=\{t\in U\,|\,s_0t=ts_0\}$.
Then $W^Us_0\subset W^{T\cup\{s_0\}}$, 
that is, $S(ws_0)=T\cup\{s_0\}$ for any $w\in W^U$.
\end{Lemma}

We first show two technical lemmas.

\begin{Lemma}\label{Lem:right-angled:2-2}
Let $(W,S)$ be an irreducible right-angled Coxeter system and let $w\in W$.
Suppose that $t_1,t_2,\dots,t_n\in S$ such that 
\begin{enumerate}
\item[(1)] $t_1\not\in S(w)$, 
\item[(2)] $o(t_it_{i+1})=\infty$ for any $i\in\{1,\dots,n-1\}$, and 
\item[(3)] $\{t_1,t_2,\dots,t_n\}=S$.
\end{enumerate}
Then $w(t_1t_2\cdots t_n)\in W^{\{t_n\}}$.
\end{Lemma}

\begin{proof}
By Lemma~\ref{Lem:right-angled:2}, 
we have that 
$$S(wt_1)=\{t\in S(w)\,|\,tt_1=t_1t\}\cup\{t_1\}=T_1\cup\{t_1\},$$
where $T_1=\{t\in S(w)\,|\,tt_1=t_1t\}$.
Also 
\begin{align*}
S(wt_1t_2)&=\{t\in S(wt_1)\,|\,tt_2=t_2t\}\cup\{t_2\} \\
&=\{t\in T_1\cup\{t_1\}\,|\,tt_2=t_2t\}\cup\{t_2\} \\
&=\{t\in T_1\,|\,tt_2=t_2t\}\cup\{t_2\} \\
&=\{t\in S(w)\,|\,tt_1=t_1t, tt_2=t_2t\}\cup\{t_2\}.
\end{align*}
By iterating the same argument, 
we obtain that 
$$S(w(t_1\dots t_n))
=\{t\in S(w)\,|\,tt_i=t_it \ \text{for any}\ i\in\{1,\dots,n\}\}\cup\{t_n\}.$$
Here $\{t_1,t_2,\dots,t_n\}=S$ by (3).
Since $(W,S)$ is irreducible, 
there does not exist $t\in S$ such that $tt_i=t_it$ 
for any $i\in\{1,\dots,n\}$.
Hence $S(w(t_1\dots t_n))=\{t_n\}$, that is, 
$w(t_1\cdots t_n)\in W^{\{t_n\}}$.
\end{proof}

\begin{Lemma}\label{lem:right-angled}
Let $(W,S)$ be an irreducible right-angled Coxeter system 
with $|\partial\Sigma(W,S)|>2$.
For each $w,v\in W$, 
there exists $x\in W$ such that 
\begin{enumerate}
\item[(1)] $\ell(x)\le 1$ and 
\item[(2)] $S(wx)\cup S(vx)\neq S$.
\end{enumerate}
\end{Lemma}

\begin{proof}
Let $w,v\in W$.
If $S(w)\cup S(v)\neq S$ then 
$x=1$ satisfies the conditions (1) and (2).
We suppose that $S(w)\cup S(v)=S$.

Then $S(w)\cap S(v)=\emptyset$.
Indeed if $s\in S(w)\cap S(v)$ then 
$st=ts$ for any $t\in S(w)$ 
by $s\in S(w)$ and 
Lemmas~\ref{lem:Coxeter-basic} and \ref{lem:right-angled:basic}~(1),
and also 
$st=ts$ for any $t\in S(v)$ 
by $s\in S(v)$ and 
Lemmas~\ref{lem:Coxeter-basic} and \ref{lem:right-angled:basic}~(1).
Hence $st=ts$ for any $t\in S(w)\cup S(v)=S$ and 
$W=W_{\{s\}}\times W_{S\setminus\{s\}}$ 
which contradicts to the assumption $(W,S)$ is irreducible.
Thus $S(w)\cap S(v)=\emptyset$.

Let $s_0\in S(w)$.
Since $(W,S)$ is an irreducible right-angled Coxeter system, 
$o(s_0t_0)=\infty$ for some $t_0\in S$ 
by Lemma~\ref{lem:right-angled:basic}~(2).
Then $t_0\in S(v)$, 
because $W_{S(w)}$ and $W_{S(v)}$ are finite 
by Lemma~\ref{lem:Coxeter-basic} 
and $W_{\{s_0,t_0\}}$ is infinite.

Here by Lemma~\ref{Lem:right-angled:2}, 
$$ S(vs_0)=\{s_0\}\cup\{t\in S(v)\,|\,ts_0=s_0t\}.$$

If $S(ws_0)\cup S(vs_0)\neq S$ then 
$x=s_0$ satisfies the conditions (1) and (2).
We suppose that $S(ws_0)\cup S(vs_0)=S$.
Since $t_0\not\in S(vs_0)$, $t_0\in S(ws_0)$.
Hence $t_0s=st_0$ for each $s\in S(ws_0)=S\setminus S(vs_0)$.
Here 
\begin{align*}
S\setminus S(vs_0)&=S\setminus(\{s_0\}\cup\{t\in S(v)\,|\,ts_0=s_0t\}) \\
&\supset S\setminus (\{s_0\}\cup S(v)) \\
&= S(w)\setminus \{s_0\}, 
\end{align*}
since $S(w)\cup S(v)=S$ and $S(w)\cap S(v)=\emptyset$.
Hence $t_0s=st_0$ for any $s\in S(w)\setminus \{s_0\}$.
Since $t_0\in S(v)$, we also have that 
$t_0s=st_0$ for any $s\in S(v)=S\setminus S(w)$.
Thus $t_0s=st_0$ for any $s\in S\setminus\{s_0\}$.
Here $t_0$ is an arbitrary element of $S$ with $o(s_0t_0)=\infty$.

Let $A=\{t\in S\,|\,o(s_0t)=\infty\}$.
Then $st=ts$ for any $t\in A$ and $s\in S\setminus\{s_0\}$ 
by the above argument.
Also $s_0t=ts_0$ for any $t\in S\setminus A$ by the definition of $A$, 
since $(W,S)$ is right-angled.
Hence we obtain that 
$$ W=W_{A\cup\{s_0\}}\times W_{S\setminus(A\cup\{s_0\})}.$$
Since $(W,S)$ is irreducible, 
$S\setminus(A\cup\{s_0\})=\emptyset$ and 
$S=A\cup\{s_0\}$.
Hence 
$$W=W_{\{s_0\}}*W_A.$$
Here $s_0$ is an arbitrary element of $S(w)$.

If there does not exist $x\in W$ which satisfies 
the conditions (1) and (2), 
then by the same argument for $t_0\in S(v)$, 
we have that 
$$W=W_{\{s_0\}}*W_{\{t_0\}}\cong \Z_2*\Z_2.$$
Then the boundary $\partial\Sigma(W,S)$ is two-points set, 
which contradicts to the assumption $|\partial\Sigma(W,S)|>2$.

Thus there exists $x\in W$ which satisfies the conditions (1) and (2).
\end{proof}

Using Theorem~\ref{Thm7-1} and lemmas above, 
we prove Theorem~\ref{Thm8-1}.

\begin{proof}[Proof of Theorem~\ref{Thm8-1}]
Let $(W,S)$ be an irreducible right-angled Coxeter system 
with $|\partial\Sigma(W,S)|>2$.
To show that the boundary $\partial\Sigma(W,S)$ is a scrambled set, 
by Theorem~\ref{Thm7-1}, we prove that 
there exist $s_0,t_0\in S$ and a number $K>0$ such that 
\begin{enumerate}
\item[(1)] $o(s_0t_0)=\infty$ and 
\item[(2)] for each $w,v\in W$, 
there exists $x\in W$ such that 
$\ell(x)\le K$ and $wx,vx\in W^{\{s_0\}}$.
\end{enumerate}

Let $s_0\in S$ and let $w,v\in W$.
By Lemma~\ref{lem:right-angled}, 
there exists $x_0\in W$ such that $\ell(x_0)\le 1$ and 
$S(wx_0)\cup S(vx_0)\neq S$.
Then there is $t_1\in S\setminus S(wx_0)\cup S(vx_0)$.
Since $(W,S)$ is an irreducible right-angled Coxeter system, 
by Lemma~\ref{lem:right-angled:basic}~(2), 
there exist $t_2,\dots,t_n\in S$ 
such that $o(t_it_{i+1})=\infty$ for each $i\in\{1,\dots,n-1\}$, 
$t_n=s_0$ and $\{t_1,t_2,\dots,t_n\}=S$.

Then by Lemma~\ref{Lem:right-angled:2-2}, 
\begin{align*}
&wx_0(t_1t_2\cdots t_n)\in W^{\{t_n\}}=W^{\{s_0\}} \ \text{and} \\
&vx_0(t_1t_2\cdots t_n)\in W^{\{t_n\}}=W^{\{s_0\}}.
\end{align*}
Hence we can take a large number $K>0$ such that 
for each $w,v\in W$, there is $x\in W$ such that 
$\ell(x)\le K$ and $wx,vx\in W^{\{s_0\}}$, 
because $S$ is finite.
Also there exists $t_0\in S$ with $o(s_0t_0)=\infty$ 
by Lemma~\ref{lem:right-angled:basic}~(2).

Therefore the boundary $\partial\Sigma(W,S)$ is a scrambled set 
by Theorem~\ref{Thm7-1}.
\end{proof}

Let $(W,S)$ be a Coxeter system. 
There exists a unique decomposition $\{S_1,\ldots,S_r\}$ of $S$ 
such that $W$ is the direct product of 
the parabolic subgroups $W_{S_1},\ldots,W_{S_r}$ and 
each Coxeter system $(W_{S_i},S_i)$ is irreducible 
(\cite{Bo}, \cite[p.30]{Hu}).
We define 
$$\tilde{S}:=\bigcup \{S_i \,|\, W_{S_i} \ \text{is infinite} \}.$$
The Coxeter system $(W,S)$ determines the subset $\tilde{S}$ of $S$.
By the definition, 
$$ W=W_{\tilde{S}}\times W_{S\setminus\tilde{S}}, $$
$W_{\tilde{S}}$ is infinite and $W_{S\setminus\tilde{S}}$ is finite.
Also it is known that 
the parabolic subgroup $W_{\tilde{S}}$ is the minimum parabolic subgroup 
of finite index in $W$ (\cite[Corollary~3.3]{Ho0}).

By Theorem~\ref{Thm8-1} and \cite[Theorem~5.1]{Ho5}, 
we obtain equivalent conditions of 
a right-angled Coxeter system whose boundary is a scrambled set.

\begin{Corollary}
Let $(W,S)$ be a right-angled Coxeter system with $|\partial\Sigma(W,S)|>2$.
Then the following statements are equivalent:
\begin{enumerate}
\item[(1)] $\partial\Sigma(W,S)$ is a scrambled set.
\item[(2)] $\partial\Sigma(W,S)$ is minimal.
\item[(3)] $(W_{\tilde{S}},\tilde{S})$ is irreducible.
\end{enumerate}
\end{Corollary}

\begin{proof}
$(2)\Leftrightarrow (3)$: 
The statements (2) and (3) are equivalent by \cite[Theorem~5.1]{Ho5}.

$(1)\Rightarrow (3)$: 
Suppose that $(W_{\tilde{S}},\tilde{S})$ is {\it not} irreducible.
Then $W_{\tilde{S}}=W_{S_1}\times W_{S_2}$ for some $S_1,S_2\subset \tilde{S}$, 
where $W_{S_1}$ and $W_{S_2}$ are infinite by the definition of $\tilde{S}$.
This implies that 
$$\Sigma(W_{\tilde{S}},\tilde{S})=
\Sigma(W_{S_1},S_1)\times\Sigma(W_{S_2},S_2)$$
and 
\begin{align*}
\Sigma(W,S)&=\Sigma(W_{\tilde{S}},\tilde{S})\times
\Sigma(W_{S\setminus\tilde{S}},S\setminus\tilde{S}) \\
&=\Sigma(W_{S_1},S_1)\times\Sigma(W_{S_2},S_2)\times
\Sigma(W_{S\setminus\tilde{S}},S\setminus\tilde{S}).
\end{align*}
Here 
$\Sigma(W_{S\setminus\tilde{S}},S\setminus\tilde{S})$ is bounded, 
because $W_{S\setminus\tilde{S}}$ is finite.
Hence 
$\Sigma(W_{S_1},S_1)\times\Sigma(W_{S_2},S_2)$ is quasi-dense in $\Sigma(W,S)$.
By Theorem~\ref{Thm:non}, 
we obtain that 
$\partial\Sigma(W,S)$ is {\it not} a scrambled set.

$(3)\Rightarrow (1)$: 
Suppose that $(W_{\tilde{S}},\tilde{S})$ is irreducible.
By Theorem~\ref{Thm8-1}, 
the boundary $\partial\Sigma(W_{\tilde{S}},\tilde{S})$ is a scrambled set.
Since 
$$\Sigma(W,S)=\Sigma(W_{\tilde{S}},\tilde{S})\times
\Sigma(W_{S\setminus\tilde{S}},S\setminus\tilde{S})$$ 
and 
$\Sigma(W_{S\setminus\tilde{S}},S\setminus\tilde{S})$ is bounded, 
$$\partial\Sigma(W,S)=\partial\Sigma(W_{\tilde{S}},\tilde{S}).$$
Therefore $\partial\Sigma(W,S)$ is a scrambled set.
\end{proof}

\section{Remark}

We can find some similarity between 
the conditions of CAT(0) groups and Coxeter groups 
whose boundaries are scrambled sets in this paper and 
the known conditions of CAT(0) groups and Coxeter groups 
whose boundaries are minimal sets 
in \cite{Ho1}, \cite{Ho2}, \cite{Ho3} and \cite{Ho5}.
The relation is unknown now in general.
The author has a question: 
Is it the case that a boundary of a CAT(0) group is a scrambled set 
if and only if it is minimal?

%

%
\end{document}